\documentclass{amsart}
\usepackage{amsmath,amssymb,amsthm,amscd,graphicx}
\usepackage{todonotes}
\usepackage{comment}
\usepackage{tikz}
\usepackage{tikz-cd}
\usepackage{caption}
\usepackage{subcaption}

\newcommand{\N}{\mathbb{N}}
\newcommand{\Z}{\mathbb{Z}}
\newcommand{\Q}{\mathbb{Q}}

\newcommand{\C}{\mathbb{C}}

\newcommand{\Pp}{\mathbb{P}}

\newcommand{\GL}{\operatorname{GL}}

\newcommand{\Aut}{\operatorname{Aut}}
\newcommand{\SL}{\operatorname{SL}}

\newcommand{\Imag}{\operatorname{Im}}
\newcommand{\Gal}{\operatorname{Gal}}

\theoremstyle{definition}
\newtheorem{definition}{Definition}

\newtheorem{remark}{Remark}
\newtheorem{example}{Example}
\newtheorem*{example*}{Example}

\theoremstyle{plain}
\newtheorem{lemma}{Lemma}
\newtheorem{theorem}{Theorem}
\newtheorem*{theorem*}{Theorem}
\newtheorem{prop}{Proposition}

\newtheorem*{question}{Question}

\usepackage{fullpage}
\usepackage{float}

\usepackage{url}
\usepackage{pdfpages}

\title{Computing the Level of a Fiber for Points on Modular Curves}
\author{Hailey Maxwell}
\address{Wake Forest University, Winston-Salem, NC 27109, USA}

\begin{document}

\maketitle

\begin{abstract}
   The modular curves in the family $X_1(N)$ for natural numbers $N$ parametrize elliptic curves over the complex numbers with a distinguished point of order $N$. The purpose of this paper is to better understand how to calculate the degrees of points on $X_1(\ell^n)$ for a prime $\ell$ and arbitrary positive integer $n$. In analogy with the definition of the level of a Galois representation, we construct a new definition: the level of a fiber of a closed point on a modular curve. Using this definition, we prove that, under certain conditions, if the degree of a point on $X_1(\ell^{k+1})$ is as large as possible given the degree of its image on $X_1(\ell^k),$ then its lifts on $X_1(\ell^n)$ have degree as large as possible for all $n > k$. We prove this result using techniques inspired by work of Lang and Trotter which gives a similar result for the image of $\ell$-adic Galois representations.  
\end{abstract}

\section{Introduction}\label{intro}
Elliptic curves are curves in projective space whose points form an additive abelian group. The modular curves in the family $X_1(N)$ for $N\in\N$ parametrize elliptic curves such that each noncuspidal point on $X_1(N)$ corresponds to an elliptic curve $E$ and a particular point $P$ on $E$ of order $N$, up to isomorphism. The purpose of this paper is to better understand how to calculate the degrees of points on $X_1(N).$

In particular, given a prime $\ell$, we want to know whether, if we have a closed point $x = [E,P] \in X_1(\ell^{k+1})$ with degree as large as possible given the degree of its image $[E,\ell P]$ on $X_1(\ell^k),$ we will have that lifts of $x$ on $X_1(\ell^n)$ are in degree as large as possible for all $n>k.$ If so, this would be an analogous feature to that of the level of the $\ell$-adic Galois representation of an elliptic curve, which, for odd primes $\ell$ and $k \geq 1$, is the first power $\ell^k$ such that the image of the mod $\ell^{k+1}$ representation is the complete preimage of the mod $\ell^k$ representation. See Lang-Trotter \cite[$\S6$]{langtrotter} and Sutherland-Zywina \cite[Lemma 3.7]{SutherlandZywina}.

\begin{figure}
    \centering
\includegraphics[width=0.5\linewidth]{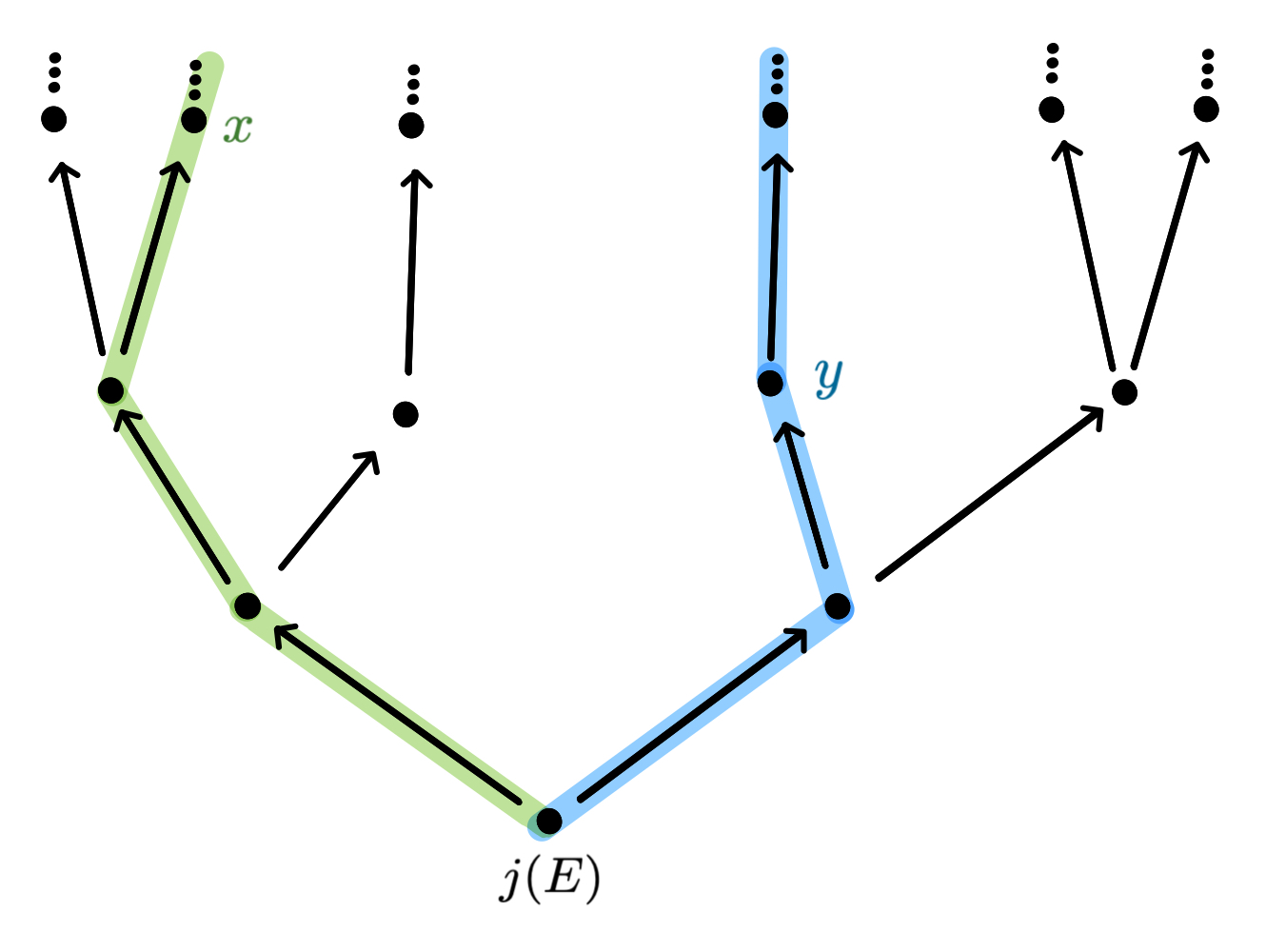}
    \caption[Some fibers of $G(E,\ell)$]{\centering Some fibers of $G(E,\ell).$ The fiber associated to $x$ is green and the fiber associated to $y$ is blue.}
    \label{introexample}
\end{figure}

Inspired by this analogy, we construct a new definition: the level of a fiber associated to a closed point on a modular curve. For a fixed prime $\ell$ and an elliptic curve $E$ we create a directed tree, called $G(E,\ell),$ using the closed points on $X_1(\ell^k)$ for all $k \geq 0$ that are lifts of the $j$-invariant of $E.$ Here, we view $j(E)$ as a closed point on $X_1(1)$ and treat it as the root of the graph. The \textbf{fibers} of $G(E,\ell)$ are the infinite paths starting at the root. Some vertices uniquely determine a fiber; that is, they and all of their descendants have degree two as vertices of the graph. For one of these vertices, the fiber through it is called its \textbf{associated} \textbf{fiber}. For an example, see Figure \ref{introexample}. As a consequence of Serre's open image theorem \cite{serre68}, if $E$ does not have complex multiplication, there are finitely many fibers in $G(E,\ell)$; see Proposition \ref{finbrancheslemma}. Finally, we define the \textbf{level of a fiber} to be the level of the point in the fiber at which the fiber stops branching. In Figure \ref{introexample}, the fiber associated to $x$ has level $\ell^3$ and the fiber associated to $y$ has level $\ell^2.$

Using this definition, we can rephrase the original question: if $x_k$ is a vertex of $G(E,\ell)$ that corresponds to a point on $X_1(\ell^k)$ and $x_k$ has only one child $x_{k+1} \in X_1(\ell^{k+1})$, does its child also have only one child? If the answer to this is yes, then by induction we have that for all $n>k,$ the vertex $x_n$ has only one child. We prove in Section \ref{definition} that this would imply that these vertices, considered as points on their corresponding modular curves, have degree growth as large as possible. In keeping with our analogy to Galois representations of elliptic curves, for $E$ defined over $\Q(j(E)),$ we prove that the level of the fiber associated to a point $[E,P] \in X_1(\ell^n)$ is at most the level of the $\ell$-adic Galois representation of $E$; see Proposition \ref{levellemma}. Then as a consequence of Serre's open image theorem \cite{serre68}, we get that for any elliptic curve $E$ without complex multiplication, the level of every fiber of $G(E,\ell)$ is finite.

Although the image of the $\ell$-adic Galois representation of $E$ can be used to impose an upper bound for the level of any fiber of $G(E,\ell)$, they need not be equal in general; see Example \ref{graphexample}. Moreover, the level of the $\ell$-adic Galois representation can be difficult to compute for elliptic curves over number fields. Thus we seek an upper bound that can be computed without needing full knowledge of the $\ell$-adic image.

Our main result answers the posed question under certain conditions:

    \begin{theorem}\label{introtheorem}
    Let $E$ be a non-CM elliptic curve over $F = \Q(j(E))$, $\ell$ be any prime, $n$ be a positive integer, and $\{P,P'\}$ be a basis for $E[\ell^{n+1}]$. Assume $P'$ is defined over $F(\ell P).$ Let $x = [E,P]$ be a closed point on $X_1(\ell^{n+1})$ such that its degree is as large as possible given the degree of its image on $X_1(\ell^n).$ Then $x$ belongs to a fiber of level at most $\ell^n.$
\end{theorem}

In Section \ref{independent_rationality}, we prove this theorem using a method inspired by Lang and Trotter \cite{langtrotter}, as appears in work of Sutherland and Zywina \cite{SutherlandZywina}. A powerful tool in the proof is the $\ell$-power map, which sends $A\to A^\ell$ for any $A \in\GL_2(\Z/\ell^n\Z)$. Roughly, when the $\ell$-power map is an injection, information at the $\ell^n$-level can be used to make conclusions at the $\ell^{n+1}$-level. However, in the context of the degrees of fibers of points on modular curves, it can be difficult to obtain such a map; see Section \ref{failed_attempt}. In particular, some additional hypotheses are needed in order for the conclusion of Theorem \ref{introtheorem} to hold; see Section \ref{counterexample}.

Using Theorem \ref{introtheorem} and a result of Bourdon, Ejder, Liu, Odumodu, and Viray \cite{BELOV}, we can also draw conclusions about isolated points on modular curves; see Proposition \ref{isolatedprop}. 

We conclude with a few examples, including an elliptic curve that satisfies the conditions of Theorem \ref{introtheorem} and an analysis of an associated fiber; see Section \ref{conditions_happen}.

\section*{Acknowledgements}

The author thanks Abbey Bourdon for her guidance throughout the research and writing processes. The author also thanks Jeremy Rouse and Hugh Howards for their comments on the thesis that formed the basis of this work. The author was partially supported by NSF grants DMS-2137659 and DMS-2145270.

\section*{Mathematica}
Some of the computations in this paper were completed using Mathematica. A file containing those computations can be found on GitHub: \\ \url{https://github.com/HMaxwell11/Computing-the-Level-of-a-Fiber-for-Points-on-Modular-Curves.git}

\section{Background} \label{background}

\subsection{Elliptic Curves}

Recall that an {elliptic curve is a nonsingular rational cubic curve in the projective plane with a distinguished rational point $\mathcal{O}.$ Unless otherwise specified, our elliptic curves are defined over an arbitrary number field $k$, or the number field $F = \Q(j(E))$. Throughout, we will refer to specific elliptic curves by their LMFDB \cite{LMFDB} labels. For a point $P = (x_0,y_0)\in E$ we define $x(P) :=x_0$ and $y(P) := y_0.$ We then define $k(P)$ to be the extension $k(x(P),y(P)).$ An elliptic curve has complex multiplication (CM) if its geometric endomorphism ring is larger than the ring of integers. 

We denote the $n$-torsion subgroup of $E$ by $E[n].$ This group is isomorphic to $\Z/n\Z\times\Z/n\Z$, so we can choose a basis $\{P,Q\}$ for the $n$-torsion points of $E.$ This basis is compatible with the multiplication-by-$m$ map between $E[n]$ and $E[\frac{n}{m}]$ for all $m \mid n$. We use $\psi_n$ to denote the $n$-th division polynomial of $E$, as implemented in Magma \cite{magma}. For $n\geq 3,$ the roots of $\psi_n$ are the $x$-coordinates of the points in $E[n]$. 

\subsection{Galois Representations}
Fix  a prime $\ell$. Given an integer $n,$ there is a natural reduction map \\ $\pi_n: \Z_\ell \to \Z/\ell^n\Z.$ For a matrix $M$, we define $\pi_n(M)$ to be the matrix whose $i,j$ entry is $\pi_n(M_{i,j}).$

\begin{definition}
    Let $E$ be an elliptic curve over $k$ and $n \in \N$. If we choose a basis $\{P,Q\}$ for $E[n] \cong \Z/n\Z\times\Z/n\Z,$ the \textbf{mod $n$ Galois representation} of $E$ is the map
    \begin{align*}\rho_{E,n}:\Gal(\overline{k}/k) &\to \Aut(E[n]) \cong \GL_2(\Z/n\Z) \\
    \sigma &\mapsto \begin{bmatrix}
        a & b \\ c & d
    \end{bmatrix}
    \end{align*}
    where $\sigma(P)=aP+cQ$ and $\sigma(Q)=bP+dQ.$
\end{definition}

If $\ell$ is a prime, the mod $\ell^n$ Galois representations of $E$ for $n \in \Z^+$ can be compiled into the $\ell$-adic Galois representation:
    $$\rho_{E,\ell^\infty}: \Gal(\overline k/k) \to \Aut(T_\ell(E))\cong \GL_2(\Z_\ell)$$
    which arises from the action of the elements of $\Gal(\overline k/k)$ on the points of $E$ with order $\ell^n$ for any $n.$

The level of the $\ell$-adic Galois representation of $E$ is the smallest integer $\ell^d$ such that $\text{im}(\rho_{E,\ell^\infty})= \pi_{d} ^{-1}(\text{im}(\rho_{E,\ell^d})).$ Due to a theorem of Serre \cite{serre68}, the level of the $\ell$-adic Galois representation of a non-CM elliptic curve is finite.

\subsection{Modular Curves}

In this section, we define the modular curves $X_1(N)$ which parametrize elliptic curves over $\C$ with points of order exactly $N$.
% \begin{theorem}\cite[Theorem 5.1 and Corollary 5.1.1.]{silvermanbook}
%     Let $E$ be an elliptic curve over $\C$. There exist a lattice $\Lambda \subset \C,$ unique up to equivalence, and a complex-analytic isomorphism $\varphi : \C/\Lambda \to E(\C)$ that respects the underlying group structures. 
% \end{theorem}
Let $\mathbb{H}$ be the upper half plane of $\C$ and $$\Gamma_1(N) := \Big\{\begin{bmatrix}
    a & b \\ c & d 
\end{bmatrix} \; | \; c \equiv 0 \pmod{N},\; a \equiv d \equiv 1 \pmod{N} \Big\} \subset \SL_2(\Z).$$ Then $Y_1(N)$ is a smooth affine curve over $\Q$ whose points are isomorphic to $\mathbb{H}/\Gamma_1(N) $ and $ X_1(N)$ is the projective closure of $Y_1(N),$ which we get by adjoining a finite number of cusps. Non-cuspidal points on $X_1(N)(\Q)$ correspond to isomorphism classes of elliptic curves over $\overline\Q$ with a point of order exactly $N.$ We write such points as $(E,P).$ We denote the associated closed point (i.e., the orbit of $(E,P)$ under the Galois action of $\overline{\Q}/\Q$) as $[E,P].$  The degree of $[E,P]$ is the size of this orbit. Alternatively, it is the degree of its residue field, which for non-CM $E$ is isomorphic to $\Q(j(E),x(P))$; see \cite[Lemma 2.5 and Remark 2.6]{BourdonNajman2021}.

There is a natural way to map a point on $X_1(N)$ to a point on $X_1(d)$ if $d\mid N $:

\begin{prop}\label{degreeprop}
For positive integers $a$ and $b$, there is a natural $\Q$-rational map $f: X_1(ab) \to X_1(a)$ defined by sending $(E,P)$ to $(E,bP).$ Moreover $$\deg(f) =  c_f\cdot b^2 \prod_{p\mid b, p\nmid a}\Bigg(1-\frac{1}{p^2}\Bigg),$$
where $c_f = \frac{1}{2}$ if $a \leq 2$ and $ab > 2$ and $c_f = 1$ otherwise.
\end{prop}

\begin{proof}
    This follows from \cite[p.66]{diamond2005}.
\end{proof}

In the case where $f$ maps $X_1(\ell^{n+1}) \to X_1(\ell^n)$ for a prime $\ell$,  this simplifies to  

   $$\deg(f) = \begin{cases}
    \ell^2 -1 \;     ,\;\;\;\ell=2 \text{ and }n = 0 ,\\
    \frac{\ell^2-1}{2}\;\;\;,\;\;\;    \ell \text { odd and } n = 0     ,\\
    \frac{\ell^2}{2}    \;\;\;\;\;\;,\;\;\;       \ell = 2 \text{ and } n = 1,\\
    \ell^2   \,\,\,\,\,\,\,\,\,\,\,  ,\;\;\;              \text{otherwise}  .
\end{cases}$$

\begin{definition}
    Let $\ell$ be a prime and let $E$ be an elliptic curve. Suppose $[E,P]\in X_1(\ell^k)$ for some $k.$ We say a closed point $[E,Q] \in X_1(\ell^{k+1})$ is a \textbf{direct lift} of $[E,P]$ if $P = \ell Q.$ We say a point $[E,R] \in X_1(\ell^n)$ for $n > k$ is a \textbf{lift} of $[E,P]$ if $P=\ell^{n-k}R.$
\end{definition}

\begin{theorem}\label{5.1,1}
    Let $E$ be a non-CM elliptic curve defined over $\Q(j(E))$ and $\ell$ a prime. Let $\ell^d$ be the level of the $\ell$-adic Galois representation of $E$ and $n \geq d$. If $x \in X_1(\ell^{n+1})$ is a closed point with $j(x) = j(E),$ then $\deg(x) = \deg(f)\deg(f(x))$ where $f$ is the map $X_1(\ell^{n+1}) \to X_1(\ell^n).$ 
\end{theorem}
\begin{proof}This is Proposition 5.8 from \cite{BELOV} restricted to the case where $f$ maps $X_1(\ell^{n+1}) \to X_1(\ell^n).$ \end{proof}

\subsection{Isolated Points}
Let $C/k$ be a curve. Let $d$ be a positive integer and denote the $d$-th symmetric power of $C$ by $C^{(d)}.$ Assume $C(k) \ne \phi$ and let $P_0\in C(k)$. Define the map
\begin{align*}
\Phi_d:C^{(d)} &\to \text{Jac}(C)\\
D &\mapsto [D-dP_0].
\end{align*}

We use this map to define what it means for a closed point $x$ of degree $d$ on $C$ to be isolated. First, we say $x$ is \textbf{$\Pp^1$-parameterized} if there exists $x'\ne x \in C^{(d)}(k)$ such that $\Phi_d(x)=\Phi_d(x').$ If no such $x'$ exists, $x$ is \textbf{$\Pp^1$-isolated}. Second, we say $x$ is \textbf{AV-parameterized} if there exists a positive rank abelian subvariety $A/k$ with $A \subset \text{Jac}(C)$ such that $\Phi_d(x)+A\subset\Imag(\Phi_d).$ If no such $x'$ exists, $x$ is \textbf{AV-isolated}. Finally, we say $x$ is \textbf{isolated} if $x$ is both \textbf{$\Pp^1$-isolated} and \textbf{AV-isolated}. There are only finitely many isolated points of any degree, and the set of such points contains $C(k)$ if genus$(C)\geq2$ \cite{faltings83}.

The following theorem provides a condition under which,  given a map between curves and an isolated point in the domain, the image of the isolated point is guaranteed to also be isolated. 

\begin{theorem}\cite[Theorem 4.3]{BELOV}\label{isolatedthm}
    Let $f:C\to D$ be a finite map of curves, let $x\in C$ be a closed point, and let $y = f(x) \in D$. Assume that $\deg(x) = \deg(y)\cdot\deg(f).$ If $x$ is isolated, then $y$ is also isolated.
\end{theorem}

\section{Connecting to the Level of a Galois Representation}

\subsection{Generalizing a Previous Result}
The following tool is a result from Sutherland and Zywina that is useful when computing Galois representations of elliptic curves. See also \cite[Proposition 3.5]{BELOV}  and \cite[page 47]{langtrotter}.

\begin{prop} \cite[Lemma 3.7]{SutherlandZywina} \label{SZprop}
Let $\ell$ be a prime and let $G$ be an open subgroup of $\GL_2(\Z_\ell).$ For each integer $m \geq 1$, let $i_m$ be the index of the image of $G$ in $\GL_2(\Z/\ell^m\Z).$ If $i_{n+1} = i_n$ for an integer $n \geq 1,$ with $n \ne 1$ if $\ell = 2,$ then $[\GL_2(\Z_\ell):G] = i_n$.
\end{prop}

This shows that the level of an $\ell$-adic Galois representation is $\ell^n$ where $n$ is the first positive integer (with $n \geq 2$ if $\ell = 2$) such that $\text{im}(\rho_{E,\ell^{n+1}})= \pi ^{-1}(\text{im}(\rho_{E,\ell^{n}})).$ Because we seek an analogous feature for points on modular curves, we generalize to the following theorem, using Remark 3.8 from the same paper \cite{SutherlandZywina}. The proof is nearly identical to the proof of Proposition \ref{SZprop}. We include a sketch for completeness and to motivate the setup of later proofs; details can be found in \cite{SutherlandZywina}.

\begin{theorem}\label{generalized}

Let $\ell$ be an odd prime. Let $A$ be a unital associative $\Z_\ell$-algebra that is torsion free and finitely generated as a $\Z_\ell$-module. Assume $A$ has $k$ generators. Let $G$ be an open subgroup of $A^\times$ and $A(\ell^m)$ be $A \pmod{\ell^m}.$ For each positive integer $m,$ let $i_m$ be the index of the image of $G$ in $A(\ell^m).$ If $i_{n+1} = i_n$ for an integer $n \geq 1,$ with $n \ne 1$ if $\ell = 2,$ then $|A^\times:G| = i_n.$
\end{theorem}

\begin{proof}
It suffices to show that $i_{m+1}=i_m$ for all $m\geq n$; we proceed by induction and the base case is given. 

Let $m \geq n$ and $G_m$ be the image of $G$ in $A(\ell^m)$. Consider the following exact sequences:

 \[
            \begin{tikzcd}        
                1 \arrow{r} & K_{m+1} \arrow{r} & A(\ell^{m+1}) \arrow{r}  & A(\ell^m) \arrow{r} & 1 \\
                1 \arrow{r} & H_{m+1} \arrow{r} \arrow[hook]{u} & G_{m+1} \arrow{r} \arrow[hook]{u} & G_m \arrow{r} \arrow[hook]{u} & 1 
            \end{tikzcd}
            \]
            
Suppose $i_{m+1} = i_m$. Then the kernels $K_{m+1}$ and $H_{m+1}$ coincide, so $|H_{m+1}| = \ell^k.$ We have that $|H_{m+2}| \leq \ell^k$ so it suffices to find an injective map from $H_{m+1}$ into $H_{m+2}.$ That map is the $\ell$-power map.
\end{proof}

\subsection{A Natural but Unsuccessful Attempt} \label{failed_attempt}

In this section, let $E$ be an elliptic curve defined over $F=\Q(j(E))$ and $P\in E(\overline F)$ have order $\ell^n$ for some postitive integer $n.$ 

The proof of Theorem \ref{generalized} involves showing that one set being of size $\ell^k$ implies another set is of the same size. For most $\ell$ and most $n>1$, the largest the degree of the extension $F(x(P))$ over $F(x(\ell P))$ can be is $\ell^2$ (see Proposition \ref{degreeprop}). To apply Theorem \ref{generalized} to make conclusions about such extensions, we need $k = 2,$ so we need $A$ to be a $\Z_\ell$-algebra that has two generators as a $\Z_\ell$-module. If we fix a basis $\{X,Y\}$ for $E[\ell^n]$ for some $n,$ then the matrices $\begin{bmatrix}
    0 & 0 \\
    1 & 0
\end{bmatrix}$ and $\begin{bmatrix}
    1 & 0 \\
    0 & 0
\end{bmatrix}$ affect $X.$ If we are only interested in what happens to the first point of our basis, a natural first attempt is to build $A$ using these two matrices. So we let

$$A = 
\Bigg\{\alpha\begin{bmatrix}
    1 & 0 \\
    0 & 0 
\end{bmatrix} + \beta\begin{bmatrix}
    0&0 \\
    1 & 0
\end{bmatrix} \;|\; \alpha,\beta \in \Z_\ell\Bigg\}.$$

While this is a natural choice given our focus, it is reasonable to suspect the argument may not work, since $A$ lacks a unity. In fact, we were able to show that the argument does not go through even for $\ell=3$.

Although the most intuitive way to apply Theorem \ref{generalized} was unsuccessful, we will see another way to use it in Section \ref{results}.

\section{The Level of a Fiber}\label{definition}

In this section we introduce the idea of the level of the fiber associated to a point on a modular curve and prove some preliminary properties. Throughout, let $E$ be an elliptic curve defined over $F = \Q(j(E)).$

\subsection{Main Construction}
\begin{definition}
    Let $\ell$ be a prime and $E$ be an elliptic curve over $F$. Let $x$ be the closed point on $X_1(1)$ that corresponds to $E.$  Let $\{x_{1,i}\}$ be the set of closed points on $X_1(\ell)$ that are lifts of $x.$ That is, each $x_{1,i}$ is $[E,P]$ for some $P\in E$ with order $\ell$. Let $\{x_{2,j}\}$ be the set of closed points on $X_1(\ell^2)$ that are lifts of any of the $x_{1,i}$. Create a directed tree with $x$ as the root and an edge pointing to each $x_{1,i}.$ For each $x_{1,i},$ add edges pointing to those $x_{2,j}$ that are direct lifts of $x_{1,i}$ (see Figure \ref{tree}). Continue: for each $n$ let $\{x_{n,i}\}$ be the set of closed points on $X_1(\ell^n)$ that are lifts of $x$ and let $\{x_{n+1,j}\}$ be the set of closed points on $X_1(\ell^{n+1})$ that are lifts of $x.$ For each $n$, create edges from $x_{n,i}$ to each $x_{n+1,j}$ that is a direct lift of $x_{n,i}.$  We will call this directed graph $G(E,\ell).$
\end{definition}

\begin{figure}[h]
    \centering    \includegraphics[width=0.55\linewidth]{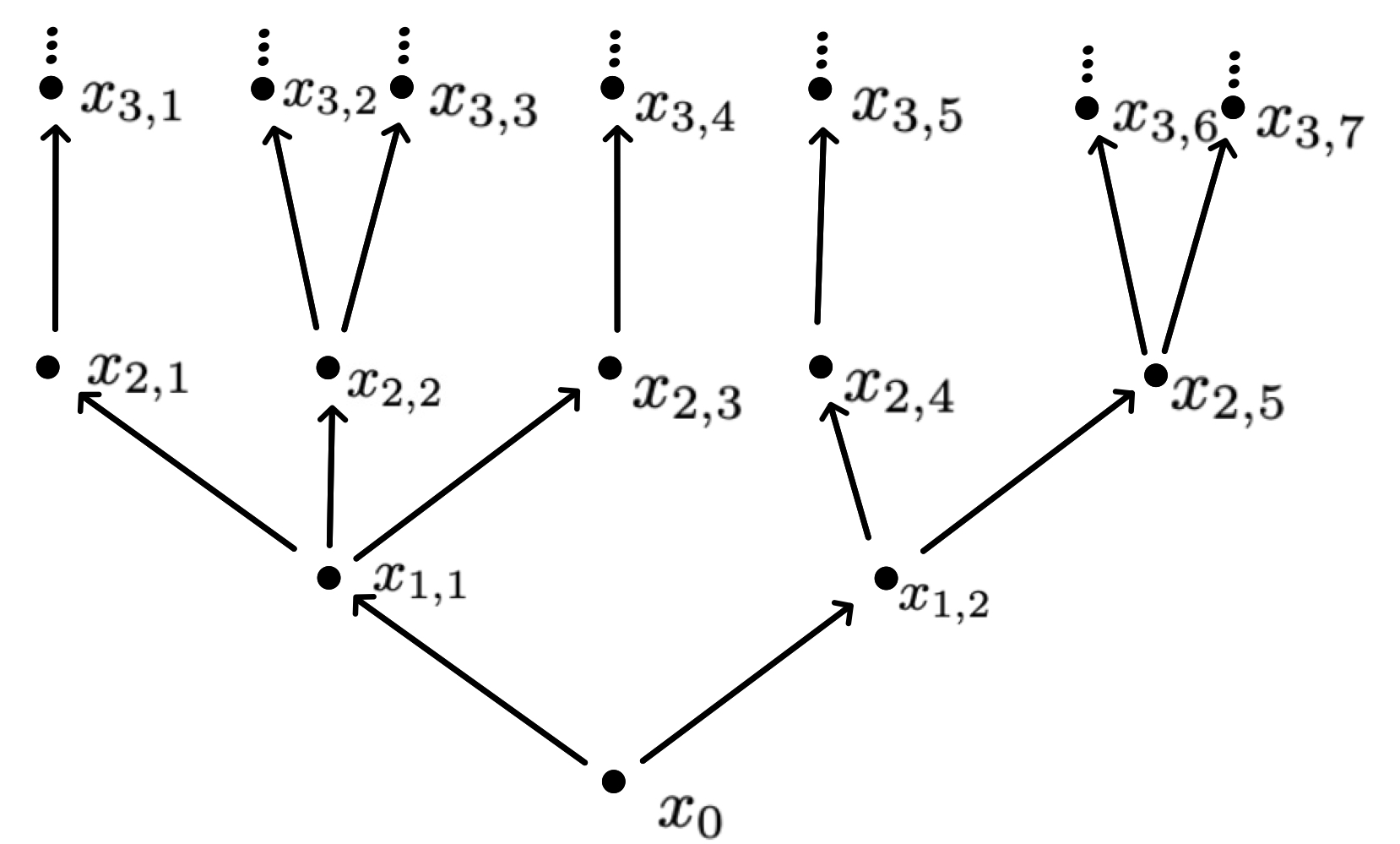}
    \caption{$G(E,\ell)$ for a generic $E$ and $\ell$}
    \label{tree}
\end{figure}

\begin{remark}
Given a point $y$ on $G(E,\ell)$ and its parent $x,$ we know that $y$ is a point on $X_1(\ell^n)$ for some $n \geq 1$ that corresponds to $[E,P]$. Then $x$ is a point on $X_1(\ell^{n-1})$ that corresponds to $[E,\ell P].$ 
\end{remark}

\begin{definition}
 Let $\ell$ be a prime, $E$ an elliptic curve, and $G = G(E,\ell).$ A  \textbf{fiber} of $G$ is an infinite path that starts at the root of $G$.
\end{definition}

\begin{definition}\label{fiberdef}
    Let $\ell$ be a prime, $E$ an elliptic curve, and $G = G(E,\ell).$ Let $S$ be the set of vertices $x \in G$ such that $x$ and all descendants of $x$ have degree 2 as vertices on $G$. For a given $x \in S$ there is a unique fiber that goes through $x$. We call that path the \textbf{fiber associated to $x$}. 
\end{definition}

\begin{example}
    In Figure \ref{fiberexample}, the green fiber is the fiber associated to $x_{3,3}$. The blue fiber is both the fiber associated to $x_{3,4}$ and the fiber associated to $x_{2,3}.$ There is no fiber associated to $x_{1,1}.$
\end{example}

\begin{figure}[h]
    \centering    \includegraphics[width=0.55\linewidth]{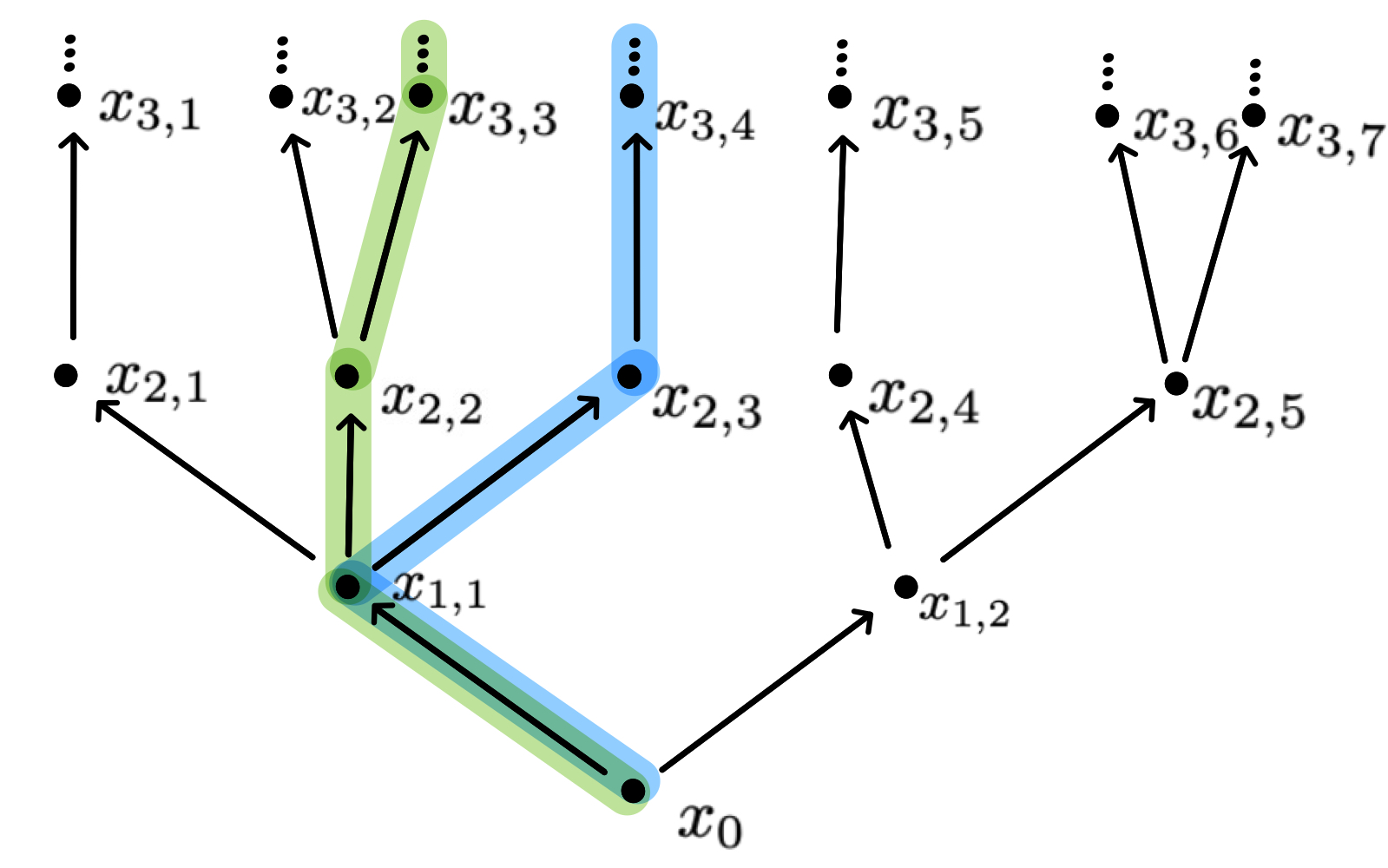}
    \caption{Some fibers of $G(E,\ell)$ for a generic $E$ and $\ell$}
    \label{fiberexample}
\end{figure}

    For a non-CM elliptic curve $E,$ this set $S$ is nonempty because, though the tree is infinite, there are only finitely many branches (see Proposition \ref{finbrancheslemma}). The uniqueness of the fiber determined by $x$ comes from the assumption that $x$ and its descendants all have degree 2. This prevents the path from branching after it reaches $x$, since the branching vertex would have to have degree at least 3.

% \begin{definition}
%     Let $S$ be defined as above. Fix $x \in S.$ Among the points in $S$ that are in the fiber associated to $x$, there is one with the shortest root-to-vertex path.  That is, one that is on $X_1(\ell^k)$ for smallest $k$ among such points. The \textbf{level of the fiber associated to $x$} is the level of that point. That is, if that point is on $X_1(\ell^k),$ then the level is $\ell^k.$
% \end{definition}

\begin{definition}
The \textbf{level of a fiber} is the level of the point in the fiber at which the fiber stops branching.
\end{definition}

\begin{example}
    In Figure \ref{fiberexample}, the green fiber has level $\ell^3$ and the blue fiber has level $\ell^2$. 
\end{example}

\begin{example}\label{graphexample}
    Let $E$ be the elliptic curve $y^2=x^3+21x+26.$ The third division polynomial of $E$ factors as two linear terms and a quadratic term. So, there are three closed points on $X_1(3)$ that correspond to $E.$ They are $x_{1,1}$ and $x_{1,2}$ with degree one and $x_{1,3}$ with degree two. 
%     The first two rows of $G(E,3)$ are depicted in Figure \ref{firsttworows}.
% 
%     \begin{figure}[h]
%     \centering    \includegraphics[width=0.35\linewidth]{pictures/firsttworows.png}
%     \caption{The bottom two rows of $G(E,3)$ for $E$: $y^2=x^3+21x+26$}
%     \label{firsttworows}
% \end{figure}
% 
    The factorization of the ninth division polynomial of $E$ over $\Q$ has eight terms. The first three are the same terms as the third division polynomial, because any point with order dividing 3 also has order dividing 9. The other five are three degree 3 terms, one degree 9 term, and one degree 18 term. So there are five closed points on $X_1(9)$ that correspond to $E$ and they have the same degrees as the corresponding terms in the factorization. The first three rows of $G(E,3)$ are depicted in Figure \ref{firstthreerows}.

    \begin{figure}[h]
    \centering    \includegraphics[width=0.55\linewidth]{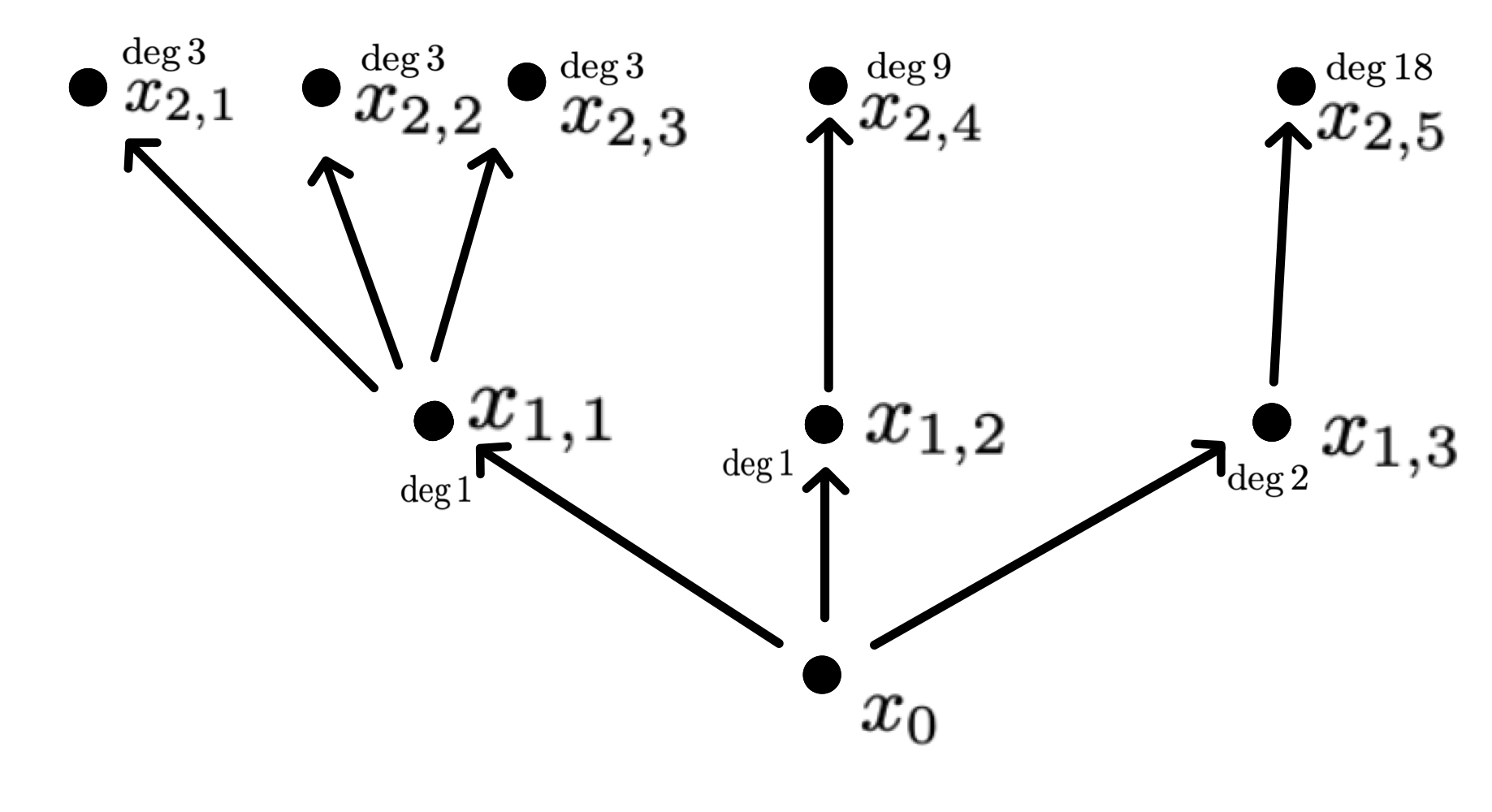}
    \caption{The bottom three rows of $G(E,3)$ for $E$: $y^2=x^3+21x+26$}
    \label{firstthreerows}
\end{figure}

    The 3-adic Galois representation of $E$ has level 9. Using Proposition \ref{levellemma}, we conclude that, for all $n > 2,$ there are five closed points on $X_1(3^n)$ that correspond to $E,$ each of which lies directly above a point on $X_1(3^{n-1}).$ Those points are $x_{n,1},...,x_{n,5}$ and their degrees can be computed using the formulas 
    \begin{align*}
    \deg(x_{n,1}) &= 3 \cdot (9)^{n-2}\\
    \deg(x_{n,2}) &= 3 \cdot (9)^{n-2}\\
    \deg(x_{n,2}) &= 3 \cdot (9)^{n-2}\\
    \deg(x_{n,2}) &= 9 \cdot (9)^{n-2}\\
    \deg(x_{n,2}) &= 18 \cdot (9)^{n-2}.
    \end{align*}
    % The degree of the map from $X_1(3^{n})\to X_1(3^{n-1})$ is 9 for all  $n > 2,$ which is where the 9 in the formulas comes from. 
    Thus Figure \ref{allrows} is the complete $G(E,3).$ There are five fibers of $G(E,3),$ here labeled $f_1,...f_5$ for ease of reference. For example, fiber $f_2$ is the purple path and fiber $f_5$ is the blue path. Fibers $f_4$ and $f_5$ have level $3$ and fibers $f_1,f_2,$ and $f_3$ have level 9.
\end{example}

\begin{figure}[h]
    \centering    \includegraphics[width=0.55\linewidth]{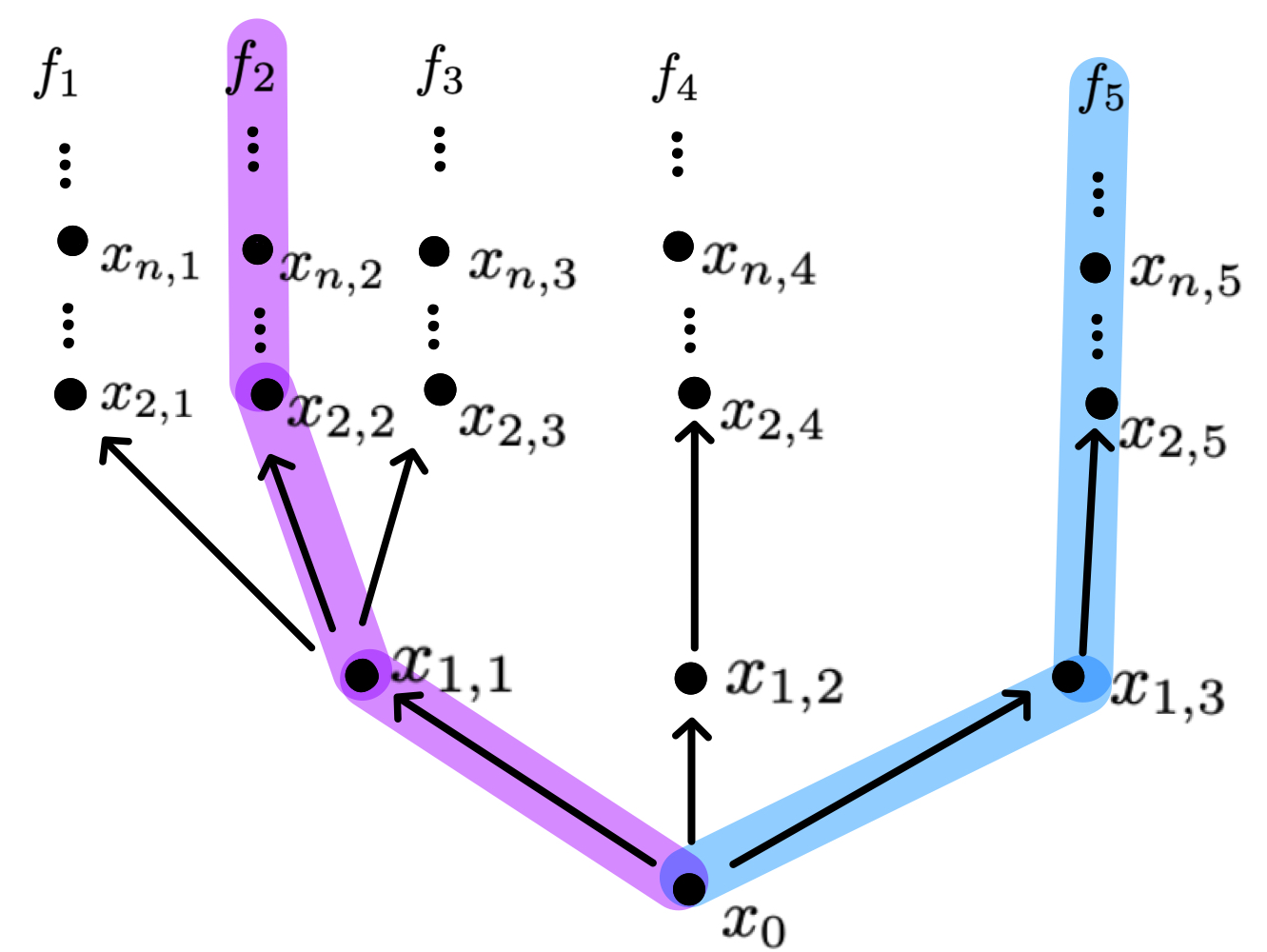}
    \caption{$G(E,3)$ for $E$: $y^2=x^3+21x+26$, with fibers identified}
    \label{allrows}
\end{figure}

\subsection{Basic Properties}
We next prove some preliminary results about this new definition. First, we show how the level of a fiber can be used to prove the existence of isolated points on some modular curves.

\begin{prop}\label{isolatedprop}
    Let $E$ be an elliptic curve over $F,$ $\ell$ a prime, and $x \in X_1(\ell^n)$ an isolated point with $j(x) = j(E).$ Suppose, as a point on $G=G(E,\ell),$ $x$ belongs to a fiber of level $\ell^k$ where $k \leq n.$ Let $f:X_1(\ell^n)\to X_1(\ell^k)$ be the natural map. Then $f(x)$ is an isolated point on $X_1(\ell^k)$.
\end{prop}

\begin{proof}
    The assumption that $x$ belongs to a fiber of level $\ell^k$ implies that $\deg(x) = \deg(f)\cdot\deg(f(x)).$ Applying Theorem \ref{isolatedthm}, we get that $x$ being isolated forces $f(x)$ to be isolated. 
\end{proof}

\begin{prop} \label{finpointsatlevel}
    Let $E$ be an elliptic curve over $F$ and $\ell$ a prime. For each $n \in \Z^+$ there are finitely many points $n$ edges away from the root of $G = G(E,\ell).$
\end{prop}

\begin{proof}
    Let $x$ be a vertex of $G$. If there are $n$ edges between $x$ and the root, then $x$ is a closed point on $X_1(\ell^n)$ that corresponds to $(E,P)$ where $P \in E[\ell^n]$.  Since $E[\ell^n] \cong \Z/\ell^n\Z \times \Z/\ell^n\Z$, there can only be finitely many such points. 
\end{proof}

\begin{prop}\label{branchlemma}
    Let $x$ be a vertex on $G(E,\ell)$ and $y$ be a child of $x.$ The degree of $y$ is as large as possible given the degree of $x$ if and only if $x$ is not a branch vertex. 
\end{prop}

\begin{proof}
    
If $x \in X_1(\ell^n)$ is a branch vertex, then it has multiple children in $ X_1(\ell^{n+1}).$ The degrees of the children sum to the degree of $x$ times the degree of the natural map $ f: X_1(\ell^{n+1}) \to X_1(\ell^{n}).$ Since each child will have degree at least the degree of $x$ (in particular, nonzero), no single child can have degree equal to the degree of $x$ times the degree of the  map, which is the largest degree a child of $x$ could possibly have. 

If $x$ is not a branch vertex, then it only has one child, $y \in X_1(\ell^{n+1})$. Then the degree of $y$ must be equal to the degree of $x$ times the degree of the map. Thus $y$ has degree as large as possible.
\end{proof}

\begin{prop}\label{levellemma}
    Let $E$ be a non-CM elliptic curve over $F$, $\ell$ a prime, and $G = G(E,\ell)$. The level of any fiber on $G$ is at most the level of the $\ell$-adic Galois representation of $E.$
\end{prop}

\begin{proof}
Let $\ell^d$ be the level of the $\ell$-adic Galois representation of $E$. Assume for the sake of contradiction that there exists a fiber of level $\ell^n$ where $n > d$. Then there exists $x \in X_1(\ell^{n-1})$ such that the point corresponding to $x$ on $G(E,\ell)$ has two children. Those children correspond to lifts of $x$ on $X_1(\ell^{n}).$ Let $y$ and $y'$ be those lifts. If $f$ is the natural projection map from $X_1(\ell^n)\to X_1(\ell^{n-1}),$ then by Theorem \ref{degreeprop}, $$\deg(y)=\deg(f)\deg(x).$$ Identically,  
     $$\deg(y')=\deg(f)\deg(x).$$ 

     Since $y$ and $y'$ are both children of $x,$  $\deg(f)\deg(x) = \deg(y)+\deg(y') = 2\deg(f)\deg(x).$ Then $\deg(f)\deg(x) = 0$, which implies $\deg(y) = 0,$ which is not possible. \qedhere

\end{proof}

\begin{remark}
    In particular, this shows that the set $S$ in Definition \ref{fiberdef} is nonempty for a non-CM elliptic curve.
\end{remark}

\begin{prop}\label{finbrancheslemma}
    Let $E$ be an elliptic curve without complex multiplication and let $\ell$ be a prime. Then there are finitely many fibers in $G = G(E,\ell).$
\end{prop}

\begin{proof}
    By Serre's open image theorem \cite{serre68}, the $\ell$-adic Galois representation of $E$ has level $\ell^d$ for some $d \in \Z^+.$  From Proposition \ref{finpointsatlevel}, there are finitely many points $d$ edges away from the root. From Proposition \ref{levellemma}, no new branches can occur more than $d$ edges away from the root of $G.$ Since there are finitely many branch vertices less than $d$ edges away from the root and no branch vertices $d$ or more edges away, there are finitely many branch vertices and thus finitely many branches. Each fiber is a branch, so there are also finitely many fibers.
\end{proof}

\section{Results} \label{results}
In this section, let $E$ be an elliptic curve defined over $F=\Q(j(E))$ and $\ell$ a prime. We use $G(E,\ell)$ to make concrete the question at the heart of this paper, rephrased below:

\begin{question}
    If $x_k$ is a vertex of $G(E,\ell)$ that corresponds to a point on $X_1(\ell^k)$ and $x_k$ has only one child, $x_{k+1} \in X_1(\ell^{k+1})$, does its child also have only one child? 
\end{question}

We partially answer this question in Proposition \ref{answer}.

In \cite{SutherlandZywina}, Sutherland and Zywina use the $\ell$-power map to draw conclusions about the level of the $\ell$-adic Galois representation of an elliptic curve. We will use the $\ell$-power map in a similar way in this section. Their proof (see Theorem \ref{generalized} for a generalization) involved using that the $\ell$-power map is an injection between two sets to show that they are the same size. In this section, we use this idea on the sets $H_{n}$ defined below. Our goal is to show that, if $H_{n+1}$ is as big as possible given the size of $H_{n},$ then $H_{n+2}$ is also as big as possible. In order to do this, we will investigate when the $\ell$-power map is an injection $H_{n+1} \to H_{n+2}.$ 

\begin{figure}[h]
    \centering
    \includegraphics[width=0.25\linewidth]{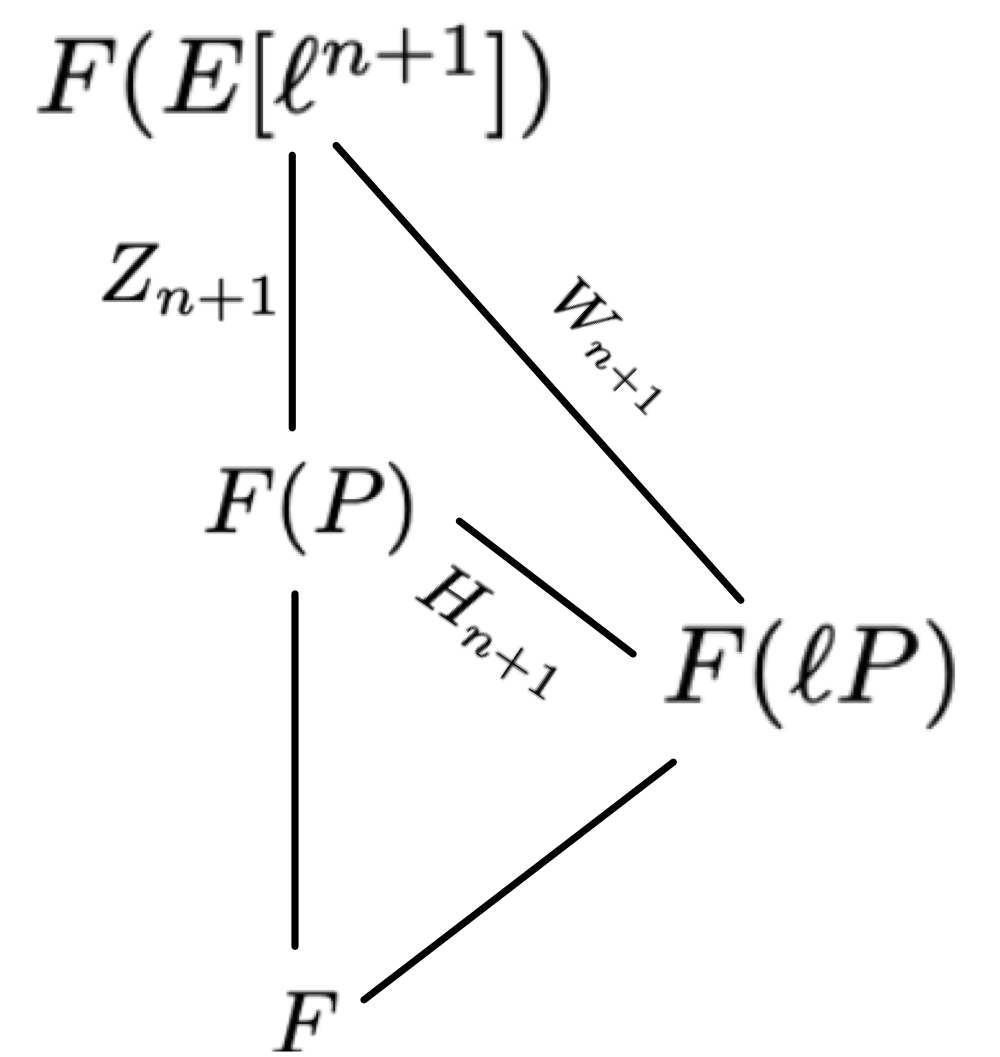}
    \caption{Field extensions for $\ell^{n+1}$}
    \label{fields}
\end{figure}

Let $E$ be an elliptic curve over $F$, $n$ be a nonnegative integer, and $\{\mathcal{P},\mathcal{Q}\}$ be a basis for $T_\ell(E)$. Then, fixing $\{P,Q\}:=\{\pi_{n+1}(\mathcal{P}),\pi_{n+1}(\mathcal{Q})\}$ as a basis for $E[\ell^{n+1}],$ we use the following notation to describe extensions of $F:$

\begin{align*}
Z_{n+1} &= \Bigg\{ \begin{bmatrix}
    1 & c \\
    0 & d
\end{bmatrix} \in \GL_2(\Z/\ell^{n+1}\Z) \Bigg\} \cap \Imag(\rho_{E,\ell^{n+1}}),\\
W_{n+1} &= \Bigg\{ \begin{bmatrix}
    1+\ell^na & c \\
    \ell^nb & 1+d
\end{bmatrix} \in \GL_2(\Z/\ell^{n+1}\Z)\Bigg\}  \cap \Imag(\rho_{E,\ell^{n+1}}),\\
H_{n+1} &= \{\alpha Z_{n+1} \;|\; \alpha \in W_{n+1}\}.
\end{align*}

The following example shows that the 3-power map is not, in general, well-defined as a map from $H_{n+1}$ to $H_{n+2}$.

\begin{example}\label{not_well_defined_example}
Let $E$ be the elliptic curve defined by $y^2=x^3-4066875x-3156806250$, with LMFDB label 50.b1 \cite{LMFDB}. Let $n = 2.$ Let $A,B \in \GL_2(\Z_3)$ such that their images in $\GL_2(\Z/27\Z)$ are 

$$A \text{ (mod 27)}= \begin{bmatrix}
    10 & 1 \\
    18 & 4
\end{bmatrix}$$
 and 
$$B \text{ (mod 27)}= \begin{bmatrix}
    10 & 6 \\
    18 & 13
\end{bmatrix}.$$
We will show that, although $A$ and $B$ represent the same cosets in $H_{n+1},$ the matrices $A^3$ and $B^3$ do not represent the same cosets in $H_{n+2}.$ 

The image of the mod $3$ Galois representation of $E$ consists of all upper triangular matrices in $\GL_2(\Z/3\Z)$. Because the $3$-adic Galois representation of $E$ has level $3,$ the mod $27$ and mod $81$ images consist of all matrices in $\GL_2(\Z/27\Z)$ and $\GL_2(\Z/81\Z),$ respectively, that reduce to an upper triangular matrix under the natural reduction maps into $\GL_2(\Z/3\Z)$. We have that $A$ (mod 3) and $B$ (mod 3) are upper triangular, so $A$ (mod 27) and $B$ (mod 27) are in the image of the mod 27 Galois representation of $E.$ 

We have  $$B^{-1}A \pmod{27}= \begin{bmatrix}
    1 & 13\\
    0 & 1
\end{bmatrix}.$$
This matrix is of the correct form to be in $Z_3$ and is in the image of the mod $27$ Galois representation of $E$ because $A$ and $B$ are, so $B^{-1}A \in Z_3.$ Thus $A$ and $B$ represent the same element of $H_3.$

We then have
$$(B^3)^{-1}A^3 \pmod{81}= \begin{bmatrix}
    28 & 48 \\
    0 & 55
\end{bmatrix}. $$

But because $28 \not\equiv 1 \pmod{81},$ $(B^3)^{-1}A^3 \pmod{81} \not\in Z_4.$ Thus $A^3$ and $B^3$ do not represent the same element of $H_4.$

\end{example}

This means the $\ell$-power map is not, in general, well-defined as a map from $H_{n+1}$ to $H_{n+2},$ but in Section \ref{independent_rationality} we see a case where it is.

\subsection{Independent Rationality}\label{independent_rationality}

We prove the following theorem in this section. 

    \begin{theorem*}%\label{introtheorem}
    Let $E$ be a non-CM elliptic curve over $F = \Q(j(E))$, $\ell$ be any prime, $n$ be a positive integer, and $\{P,P'\}$ be a basis for $E[\ell^{n+1}]$. Assume $P'$ is defined over $F(\ell P).$ Let $x = [E,P]$ be a closed point on $X_1(\ell^{n+1})$ such that its degree is as large as possible given the degree of its image on $X_1(\ell^n).$ Then $x$ belongs to a fiber of level at most $\ell^n.$
\end{theorem*}

 Throughout this section let $E$ be an elliptic curve and let $\{P,P'\}$ be a basis for $E[\ell^{n+1}].$ We use the notation for $H_i,$ $W_i,$ and $Z_i$ from the beginning of Section \ref{results}. Assume $P'$ is defined over $F(\ell P).$ Note that this means $F(E[\ell^{n+1}]) = F(P),$ so $$Z_{n+1} = \Bigg\{\begin{bmatrix}
    1 & 0 \\ 0 & 1
\end{bmatrix}\Bigg\} \subset \GL_2(\Z/\ell^{n+1}\Z)$$. 

\begin{lemma}
    Assuming $P'$ is defined over $F(\ell P),$ the $\ell$-power map from $H_{n+1}$ to $H_{n+2}$ is well-defined. 
\end{lemma}

\begin{proof}
    Let $A,B\in \GL_2(\Z_\ell)$ such that their images $\GL_2(\Z/\ell^{n+1}\Z)$ represent the same element of $H_{n+1}.$ That is, $\pi_{n+1}(A),\pi_{n+1}(B) \in W_{n+1}$ and $\pi_{n+1}(A)Z_{n+1} = \pi_{n+1}(B)Z_{n+1}.$ Then 
    \begin{align*}
         \pi_{n+1}(A)=\pi_{n+1}(B) &\implies \pi_{n+1} 
         (A^\ell)=\pi_{n+1}(B^\ell)\\
        &\implies \pi_{n+2}(A)=\pi_{n+2}(B)\\
        &\implies \pi_{n+2}(A)Z_{n+2}=\pi_{n+2}(B)Z_{n+2}.
    \end{align*}
   So $\pi_{n+2}(A)$ and $\pi_{n+2}(B)$ represent the same element of $H_{n+2}$, so the $\ell$-power map is well-defined.
\end{proof}

\begin{lemma}\label{raisingA}
    Let $n$ be a positive integer and $A \in \GL_2(\Z_\ell)$ such that the image of $A$ in $\GL_2(\Z/\ell^{n+1}\Z)$ is $\begin{bmatrix}
    1 + \ell^n a & 0 \\
    \ell^n b & 1 
    \end{bmatrix}.$ Then the image of $A^\ell$ in $\GL_2(\Z/\ell^{n+2}\Z)$ is $\begin{bmatrix}
    1 + \ell^{n+1} a & 0 \\
    \ell^{n+1} b & 1 
\end{bmatrix}.$
\end{lemma}

\begin{proof}
    Because every matrix in this lemma is being considered in either $\GL_2(\Z/\ell^{n+1}\Z)$ or $\GL_2(\Z/\ell^{n+2}\Z),$ we reduce them mod $\ell^{n+2}$ after every operation.
    
    Write $A$ as its image in $\GL_2(\Z/\ell^{n+2}\Z):$ 
    $$A\text{ (mod } \ell^{n+2}) = \begin{bmatrix}
    1 + \ell^n a & \ell^{n+1} j \\
    \ell^n b & 1 + \ell^{n+1}  k
\end{bmatrix}$$ for some $j,k \in \Z. $ 

We will show that, for all $w \in \Z^+, $

$$A^w \text{ (mod } \ell^{n+2})= \begin{bmatrix}
    1 + w\ell^n a + \binom{w}{2}\ell^{2n}a^2 & w\ell^{n+1} j \\
    w\ell^n b +\binom{w}{2}\ell^{2n}ab& 1 + w\ell^{n+1}  k
\end{bmatrix}.$$

% Note that $$ \begin{bmatrix}
%     1 + \ell\cdot\ell^n a + \binom{\ell}{2}\ell^{2n}a^2 & \ell\cdot\ell^{n+1} j \\
%     \ell\cdot\ell^n b +\binom{w}{2}\ell^{2n}ab& 1 + \ell\cdot\ell^{n+1}  k
% \end{bmatrix} \equiv \begin{bmatrix}
%     1 +\ell^{n+1} a & 0 \\
%     \ell^{n+1} b & 1
% \end{bmatrix} \pmod{\ell^{n+2}}$$
% because $\ell \mid \binom{\ell}{2}.$ So it suffices to show that $$A^\ell \text{ (mod } \ell^{n+2})=  \begin{bmatrix}
%     1 + \ell\cdot\ell^n a + \binom{\ell}{2}\ell^{2n}a^2 & \ell\cdot\ell^{n+1} j \\
%     \ell\cdot\ell^n b +\binom{\ell}{2}\ell^{2n}ab& 1 + \ell\cdot\ell^{n+1}  k
% \end{bmatrix}$$
% for all primes $\ell.$ 

We proceed by induction. Consider the base case:

\begin{align*}
    A^2 \text{ (mod } \ell^{n+2})&= \begin{bmatrix}
    1 + 2\ell^n a + \ell^{2n}a^2 & 2\ell^{n+1} j \\
    2\ell^n b +\ell^{2n}ab& 1 + 2\ell^{n+1}  k
\end{bmatrix}\\
&= \begin{bmatrix}
    1 + 2\ell^n a + \binom{2}{2}\ell^{2n}a^2 & 2\ell^{n+1} j \\
    2\ell^n b +\binom{2}{2}\ell^{2n}ab& 1 + 2\ell^{n+1}  k
\end{bmatrix}.
\end{align*}

Now assume the statement holds for $w.$ Then 
\begin{align*}
A^{w+1} \pmod{\ell^{n+2}}%&= A\cdot A^w\pmod{3^{n+2}}\\
&= \begin{bmatrix}
    1 + \ell^n a & \ell^{n+1} j \\
    \ell^n b & 1 + \ell^{n+1}  k
\end{bmatrix}
\begin{bmatrix}
    1 + w\ell^n a + \binom{w}{2}\ell^{2n}a^2 & w\ell^{n+1} j \\
    w\ell^n b +\binom{w}{2}\ell^{2n}ab& 1 + w\ell^{n+1}  k
\end{bmatrix} \\
% & = \begin{bmatrix}
%     1 + (w+1)\ell^n a +(w+ \binom{w}{2})\ell^{2n}a^2 &(w+1)\ell^{n+1} j \\
%    (w+1)\ell^n b +(w+\binom{w}{2})\ell^{2n}ab& 1 + (w+1)\ell^{n+1}  k
% \end{bmatrix} \\
& = \begin{bmatrix}
    1 + (w+1)\ell^n a +\binom{w+1}{2}\ell^{2n}a^2 &(w+1)\ell^{n+1} j \\
   (w+1)\ell^n b +\binom{w+1}{2}\ell^{2n}ab& 1 + (w+1)\ell^{n+1}  k
\end{bmatrix}.
\end{align*}

Since the statement holds for all positive integers, it holds for $\ell,$ as desired.

Finally, note that, for odd $\ell$, because $\ell \mid \binom{\ell}{2}$, 
$$ A^{\ell} \equiv \begin{bmatrix}
    1 + \ell\cdot\ell^n a + \binom{\ell}{2}\ell^{2n}a^2 & \ell\cdot\ell^{n+1} j \\
    \ell\cdot\ell^n b +\binom{\ell}{2}\ell^{2n}ab& 1 + \ell\cdot\ell^{n+1}  k
\end{bmatrix} \equiv \begin{bmatrix}
    1 +\ell^{n+1} a & 0 \\
    \ell^{n+1} b & 1
\end{bmatrix} \pmod{\ell^{n+2}}.\qedhere$$
\end{proof}

\begin{lemma} \label{isinjection}
    Assuming $P'$ is defined over $F(\ell P),$ the $\ell$-power map is an injection from $H_{n+1}$ into $H_{n+2}.$
\end{lemma}

\begin{proof}
    Let $A,B \in \GL_2(\Z_\ell)$ such that $A \pmod{\ell^{n+1}}$ and $B \pmod{\ell^{n+1}}$ represent elements in $H_{n+1}$ and $A^\ell \pmod{\ell^{n+2}}$ and $B^\ell \pmod{\ell^{n+2}}$ represent the same element of $H_{n+2}.$ That is, $\big(\pi_{n+2}(B^\ell)\big)^{-1}\pi_{n+2}(A^\ell) \in Z_{n+2}.$

    Because $A$ and $B$ are being considered only based on their images in $\GL_2(\Z/\ell^{n+1}\Z)$ and $\GL_2(\Z/\ell^{n+2}\Z),$ we reduce mod $\ell^{n+2}$ throughout.

    Write  
$$A \text{ (mod } \ell^{n+2})=  \begin{bmatrix}
    1 + \ell^n a & \ell^{n+1} j \\
    \ell^n b & 1 + \ell^{n+1}  k
\end{bmatrix}$$ for some $j, k \in \Z$ and $$B \text{ (mod } \ell^{n+2})= \begin{bmatrix}
    1 + \ell^n x & \ell^{n+1} z \\
    \ell^n y & 1 + \ell^{n+1} w
\end{bmatrix}$$ for some $z, w \in \Z$.

Using Lemma \ref{raisingA} and Mathematica, we compute 

% have
%     $$A^\ell \text{ (mod } \ell^{n+2})= \begin{bmatrix}
%         1+\ell^{n+1}a & 0 \\
%         \ell^{n+1}b & 1
%     \end{bmatrix}$$

%     and
%     $$B^\ell \text{ (mod } \ell^{n+2})= \begin{bmatrix}
%         1+\ell^{n+1}x & 0 \\
%         \ell^{n+1}y & 1
%     \end{bmatrix}.$$

% The inverse of  $B^\ell$ in $\GL_2(\Z/\ell^{n+2}\Z) $ is

% $$(B^\ell)^{-1} \text{ (mod } \ell^{n+2})= \begin{bmatrix}
%        (1+\ell^{n+1}x)^{-1} & 0\\
%        -\ell^{n+1}y(1+\ell^{n+1}x)^{-1} & 1
        
%     \end{bmatrix}.$$
% %verified in mathematica file "general ell"

% So 

$$
(B^\ell)^{-1}A^\ell \text{ (mod } \ell^{n+2})= \begin{bmatrix}
        (1+\ell^{n+1}a)(1+\ell^{n+1}x)^{-1} & 0 \\
        \ell^{n+1}b-\ell^{n+1}y(1+\ell^{n+1}a)(1+\ell^{n+1}x)^{-1} & 1
    \end{bmatrix}.
$$
%verified in mathematica file "general ell"

Our assumption is that

$$(B^\ell)^{-1}A^\ell\text{ (mod } \ell^{n+2})\in Z_{n+2} \subseteq \Big\{ \begin{bmatrix}
    1 & c \\
    0 & d
\end{bmatrix} |\; c,d \in \Z/\ell^{n+2}\Z \Big\},$$ 

which forces $1+\ell^{n+1}a \equiv 1+\ell^{n+1}x $ and $\ell^{n+1}b(1+\ell^{n+1}x) \equiv \ell^{n+1}y(1+\ell^{n+1}a)\pmod{\ell^{n+2}}.$ From the first congruence, we get that $a \equiv x \pmod{\ell}.$ From the second we get that $b \equiv y \pmod{\ell}.$

Since $a,b,x,$ and $y$ are defined mod $\ell,$ and $j,k,z$ and $w$ disappear when considering the images of $A$ and $B$ modulo $\ell^{n+1},$ we have that $\pi_{n+1}(A)$ and $\pi_{n+1}(B)$ represent the same element of $H_{n+1}.$
\end{proof}

\begin{prop}\label{n+1 -> n+2}
     Assume $P'$ is defined over $F(\ell P).$ If $H_{n+1}$ is as big as possible considering the size of $H_{n},$ then $H_{n+2}$ is as big as possible considering the size of $H_{n+1}.$
\end{prop}

\begin{proof}
    From the assumption that $H_{n+1}$ is as big as possible, we have that $|H_{n+1}| = \ell^2.$ Because the $\ell$-power map is an injection from $H_{n+1}$ to $H_{n+2}$ (see Lemma \ref{isinjection}), we have $|H_{n+1}| \leq |H_{n+2}|.$ Thus $|H_{n+2}|\geq \ell^2.$ Since $\ell^2$ is the biggest that $|H_{n+2}|$ can be, we have that $|H_{n+2}|$ is as big as possible.
\end{proof}

\begin{prop} \label{degree of lift}
     Let $E$ be a non-CM elliptic curve over $F$, $\ell$ be any prime, and $\{P,P'\}$ be a basis for $E[\ell^{n+1}]$. Assume $P'$ is defined over $F(\ell P).$ Let $x = [E,P]\in X_1(\ell^{n+1})$ and let $[E,Q] \in X_1(\ell^{n+2})$ be a lift of $[E,P].$ If the degree of $x$ is as large as possible given the degree of its image on $X_1(\ell^n),$ then the degree of $[E,Q] \in X_1(\ell^{n+2})$ is as big as possible. 
\end{prop}

\begin{proof}
     Let $x = [E,P]\in X_1(\ell^{n+1})$ and $[E,Q] \in X_1(\ell^{n+2})$ be a lift of $[E,P].$ Assume the degree of $x$ is as large as possible given the degree of its image on $X_1(\ell^n).$ Then the degree of the extension $F(P)$ over $F(\ell P)$ is as big as possible, so $H_{n+1}$ is as big as possible. From Proposition \ref{n+1 -> n+2}, we have that $H_{n+2}$ is as big as possible, so the degree of the extension $F(Q)$ over $F(P)$ is as big as possible, so the degree of $[E,Q] \in X_1(\ell^{n+2})$ is as big as possible, as desired. 
\end{proof}

\begin{prop}\label{answer}
     Let $E$ be a non-CM elliptic curve over $F$, $\ell$ be any prime, and $\{P,P'\}$ be a basis for $E[\ell^{n+1}]$. Assume $P'$ is defined over $F(\ell P).$ Let $x$ be the vertex in $G(E,\ell)$ that corresponds to  $[E,P]\in X_1(\ell^{n+1}).$ Let $[E,Q]$ be a lift of $[E,P],$ and let $y$ be the vertex that corresponds to $[E,Q] \in X_1(\ell^{n+2})$. If $x$ is the only child of its parent, then $y$ is the only child of $x$.
\end{prop}

\begin{proof}
    Let $E$ be a non-CM elliptic curve over $F$, $\ell$ be any prime, and $\{P,P'\}$ be a basis for $E[\ell^{n+1}]$. Assume $P'$ is defined over $F(\ell P).$ Let $x$ be the vertex in $G(E,\ell)$ that corresponds to  $[E,P]\in X_1(\ell^{n+1}).$ Let $[E,Q]$ be a lift of $[E,P],$ and let $y$ be the vertex that corresponds to $[E,Q] \in X_1(\ell^{n+2})$. Assume $x$ is the only child of its parent. Then by Proposition \ref{branchlemma}, $[E,P]$ has degree as large as possible. Then by Lemma \ref{degree of lift}, $[E,Q]$ also has degree as large as possible. Applying Proposition \ref{branchlemma} again, we get that $x$ is not a branch vertex, thus has only one child, $y$. 
\end{proof}

We are now ready to prove the theorem stated at the beginning of this section:

    \begin{theorem}\label{maintheorem}
    Let $E$ be a non-CM elliptic curve over $F = \Q(j(E))$, $\ell$ be any prime, $n$ be an integer at least 1, and $\{P,P'\}$ be a basis for $E[\ell^{n+1}]$. Assume $P'$ is defined over $F(\ell P).$ Let $x = [E,P]$ be a closed point on $X_1(\ell^{n+1})$ such that its degree is as large as possible given the degree of its image on $X_1(\ell^n).$ Then $x$ belongs to a fiber of level at most $\ell^n.$
\end{theorem}
\begin{proof}
    Let $x = [E,P]$ be a closed point on $X_1(\ell^{n+1})$ such that its degree is as large as possible given the degree of its image on $ X_1(\ell^n).$ Let $k$ be the level of the fiber associated to $x$.

    Suppose for sake of contradiction that $k > \ell^n.$ Then there exists $b>n$ such that $x_b \in X_1(\ell^b)$ is a lift of $x_{b-1} \in X_1(\ell^{b-1}),$ $x_b$ does not have degree as large as possible given the degree of $x_{b-1},$ and both $x_b$ and $x_{b-1}$ are lifts of $x.$ Let $b$ be the smallest such number that is greater than $n.$  

   Then $b-1$ is not such a number, so every point on $X_1(\ell^{b-1})$ that is a lift of $x$ must have degree as large as possible given the degree of its image on $X_1(\ell^{b-2}).$ Thus $x_{b-1}$ has degree as large as possible. Then by Proposition \ref{degree of lift}, we have that the degree of $x_b$ must also be as big as possible, contradicting that $b$ is a number of the form described. Thus no such number exists, so the level of the fiber associated to $x$ is at most $\ell^n.$
\end{proof}

\subsection{A Curve Satisfying the Conditions of Theorem \ref{maintheorem}} \label{conditions_happen}

Consider the curve $E: y^2+xy+y=x^3-x^2-14x+29$ with LMFDB label 54.b2 \cite{LMFDB}. Then $j(E) = -1167051/512,$ so $F=\Q(j(E)) = \Q.$ This curve has a point, $P'$, of order 9 defined over $F$. Let $\{P,P'\}$ be a basis for $E[9].$ Then $3 P$ is a point of order 3, and we have that $P'$ is defined over $F(3 P),$ because it is already defined over $F.$

We also need that $[E,P]$ has degree as large as possible given the degree of $[E,3P]$. That is, the extension $F(x(P))/$$F(x(3P))$ has degree $3^2 = 9.$ The factorization over $\Q$ of the third division polynomial of $E$ is $(x-1)(x^3-27x+90).$ This means that on $X_1(3),$ there are two closed points that correspond to $E$, one with degree 1 and one with degree 3. When the ninth division polynomial is divided by the third, the resulting polynomial factors over $\Q$ as three degree 1 terms, one degree 6 term, and one degree 27 term. So there are five closed points on $X_1(9)$ that correspond to $E,$ one of which has degree 27. Since the natural map between $X_1(3)$ and $X_1(9)$ has degree $3^2=9$, the only way to have a degree 27 point on $X_1(9)$ is to have it be a lift of the degree 3 point, which was $[E,3P]$. So, $[F(x(P)):F(x(3P))] = 9$. 

% Let $a = [F(P):F(x(P))]$ and $b = [F(P):F(3P)]$ (see Figure \ref{ellis3example}). This leads to the equation $27a=6b \implies 9a = 2b \implies 2 | a \text{ and } 9 | b.$ From Remark \ref{degtworemark}, we know that $a$ is either 1 or 2, so we must have $a=2$. Then we have the equation $54=6b,$ so $b = 9,$ and the second condition is satisfied.

% \begin{figure}[h]

%     \centering
%     \includegraphics[width=0.40\linewidth]{pictures/field-diagram-for-ell-is-3-example.jpeg}
%     \caption{Field extensions for $E: y^2+xy+y=x^3-x^2-14x+29$}
%     \label{ellis3example}
% \end{figure}

Note that the degree-27 point is on $X_1(9)$ and has degree as large as possible given its image on $X_1(3),$ which must be the degree-3 point. Then Theorem \ref{maintheorem} says that the degree-27 point belongs to a fiber of level at most 3. Since the degree-3 point does not have degree as large as possible on $X_1(3),$ the level is at least 3, thus the level of the fiber containing the degree-27 point is 3 (see Figure \ref{examplefibers}).

\begin{figure}[h]
    \centering
    \includegraphics[width=0.5\linewidth]{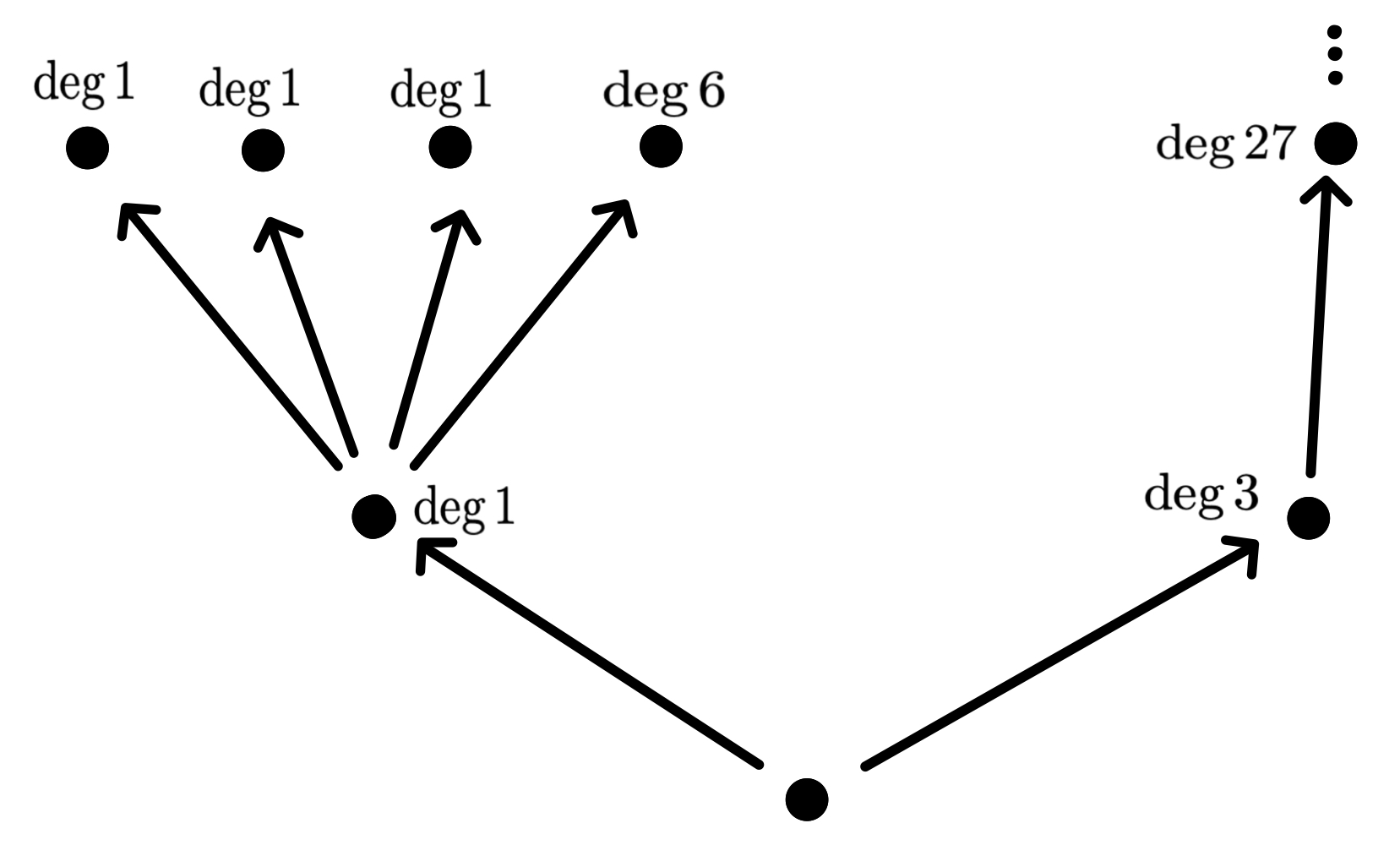}
    \caption{Some of $G(E,3)$}
    \label{examplefibers}
\end{figure}

\subsection{A Counterexample to Generality}\label{counterexample}
Let $E$ be the elliptic curve $y^2=x^3-2x+1.$ We have that $E$ does not have complex multiplication and has $2$-adic level 16 \cite{LMFDB}. Using Magma, we compute the second, fourth, eighth, and sixteenth division polynomials of $E$ to construct the first five rows of $G(E,2)$ (see Figure \ref{counterexamplepic}). Note that the vertex $x$ representing a closed point on $X_1(4)$ has only only child, $y,$ which represents a closed point on $X_1(8)$ with degree as large as possible, but $y$ has two children, $a$ and $b$ which each represent a point on $X_1(16)$ with degree smaller than the largest possible. Since $E$ has no complex multiplication, we can apply Proposition \ref{levellemma} to conclude that the fiber associated to $a$ and the fiber associated to $b$ both have level 16. 

\begin{figure}[h!]
    \centering    \includegraphics[width=0.5\linewidth]{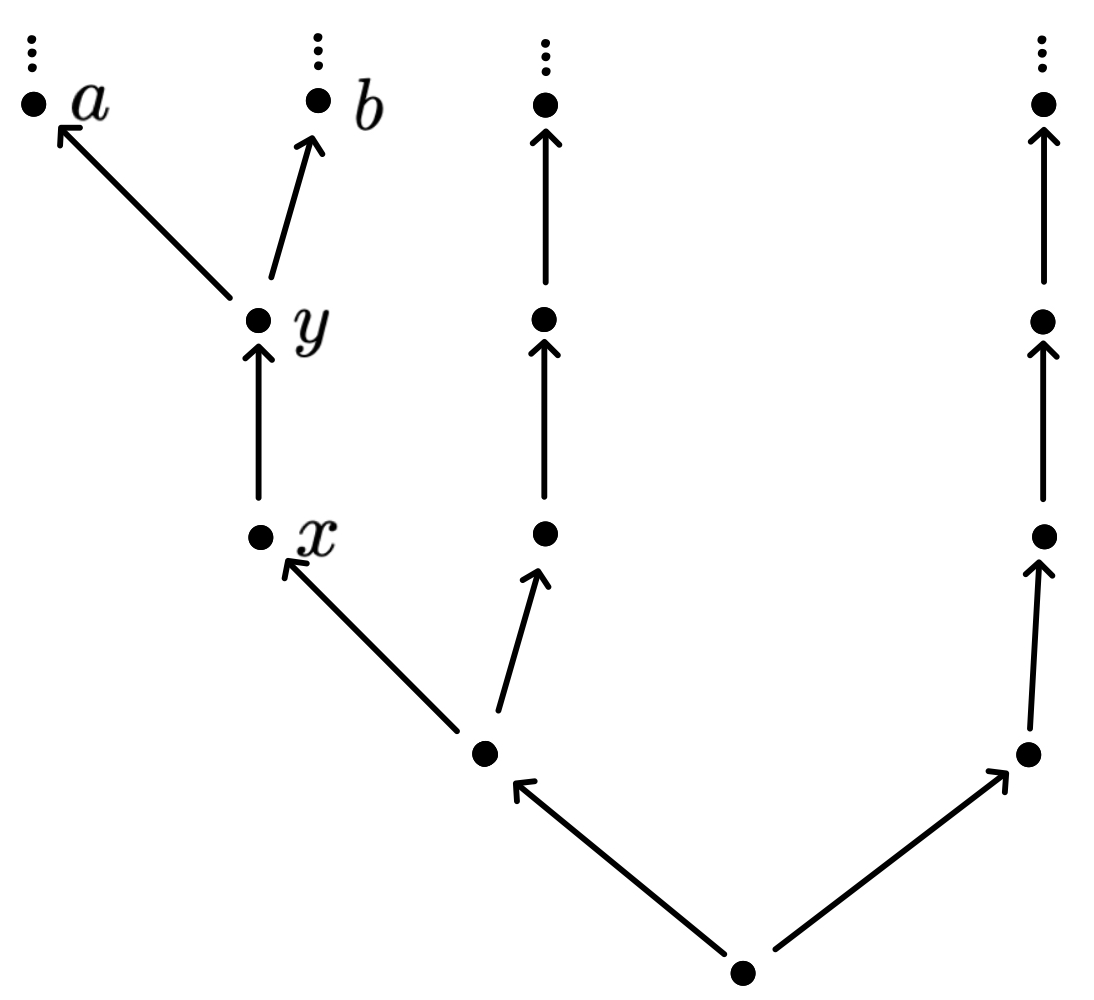}
    \caption{A counterexample to the general case}
    \label{counterexamplepic}
\end{figure}

This counterexample shows that the conclusion of Theorem \ref{maintheorem} does not hold in complete generality. Thus, some conditions are necessary. Whether the conditions stated in Theorem \ref{maintheorem} are necessary or merely sufficient is the subject of future work.

% \subsection{Case where the argument doesn't work}\label{argument_doesn't_work}
% Let $E$ be an elliptic curve such that the level of the 3-adic image of the Galois representation of $E$ is 3. Then for every $M \in \GL_2(\Z/3\Z)$ there is an $N \in \text{im}(\rho_{E,3^\infty})$ such that $N$ reduces to $M$ under the natural projection map. 

\bibliographystyle{amsplain}
\bibliography{bibfile}

\providecommand{\bysame}{\leavevmode\hbox to3em{\hrulefill}\thinspace}
\providecommand{\MR}{\relax\ifhmode\unskip\space\fi MR }
\providecommand{\MRhref}[2]{%
  \href{http://www.ams.org/mathscinet-getitem?mr=#1}{#2}
}
\providecommand{\href}[2]{#2}
\begin{thebibliography}{1}

\bibitem{magma}
Wieb Bosma, John Cannon, and Catherine Playoust, \emph{The {M}agma algebra system. {I}. {T}he user language}, J. Symbolic Comput. \textbf{24} (1997), no.~3-4, 235--265, Computational algebra and number theory (London, 1993). \MR{MR1484478}

\bibitem{BELOV}
Abbey Bourdon, \"{O}zlem Ejder, Yuan Liu, Frances Odumodu, and Bianca Viray, \emph{On the level of modular curves that give rise to isolated {$j$}-invariants}, Adv. Math. \textbf{357} (2019), 106824, 33.

\bibitem{BourdonNajman2021}
Abbey Bourdon and Filip Najman, \emph{Sporadic points of odd degree on {$X_1(N)$} coming from $\mathbb{Q}$-curves}, Preprint available at arxiv.org:2107.10909.

\bibitem{diamond2005}
F.~Diamond and J.~Shurman, \emph{A first course in modular forms}, Graduate Texts in Mathematics, Springer, 2005.

\bibitem{faltings83}
G.~Faltings, \emph{Endlichkeitss\"{a}tze f\"{u}r abelsche {V}ariet\"{a}ten \"{u}ber {Z}ahlk\"{o}rpern}, Invent. Math. \textbf{73} (1983), no.~3, 349--366. \MR{718935}

\bibitem{langtrotter}
S.~Lang and H.~Trotter, \emph{Frobenius distributions in gl2-extensions: Distribution of frobenius automorphisms in gl2-extensions of the rational numbers}, Lecture Notes in Mathematics, Springer, 1976.

\bibitem{LMFDB}
The {LMFDB Collaboration}, \emph{The {L}-functions and modular forms database}, \url{https://www.lmfdb.org}, 2023, [Online; accessed August 2023].

\bibitem{serre68}
J.P. Serre, \emph{Abelian $l$-adic representations and elliptic curves: Mcgill university lecture notes}, Mathematics lecture note series, W. A. Benjamin, 1968.

\bibitem{SutherlandZywina}
Andrew~V. Sutherland and David Zywina, \emph{Modular curves of prime-power level with infinitely many rational points}, Algebra Number Theory \textbf{11} (2017), no.~5, 1199--1229. \MR{3671434}

\end{thebibliography}

% \addcontentsline{toc}{chapter}{Curriculum Vitae}
% \includepdf[pages=-,pagecommand={}]{thesis-cv.pdf}

\end{document}